\documentclass[11pt,amsfonts]{article}
\usepackage{graphicx}
\usepackage{latexsym}
\usepackage{amssymb}
\usepackage{amsmath}
\usepackage{layout}
\newtheorem{prop}{Proposition}
\newtheorem{lemma}{Lemma}

\newtheorem{corollary}{Corollary}

\newtheorem{remark}{Remark}

\def\real{{\mathord{{\rm I\kern-2.8pt R}}}}        
\def\inte{{\mathord{{\rm I\kern-2.8pt N}}}}

\def\sZZ{{\rm Z\kern-2.8ptem{}Z}}

\def\z{{\mathchoice
  {\sZZ}
  {\sZZ}
  {\rm Z\kern-0.30em{}Z}
  {\rm Z\kern-0.25em{}Z} }}
\def\sQQ{{\kern 0.27em \vrule height1.45ex width0.03em depth0em
          \kern-0.30em \rm Q}}
\def\qu{{\mathchoice
    {\sQQ}
    {\sQQ}
  {\kern 0.225em \vrule height1.05ex width0.025em depth0em \kern-0.25em \rm Q}
  {\kern 0.180em \vrule height0.78ex width0.020em depth0em \kern-0.20em \rm Q}
        }}
\def\sCC{{\kern 0.27em \vrule height1.45ex width0.03em depth0em
          \kern-0.30em \rm C}}
\def\complex{{\mathchoice
    {\sCC}
    {\sCC}
  {\kern 0.225em \vrule height1.05ex width0.025em depth0em \kern-0.25em \rm C}
  {\kern 0.180em \vrule height0.78ex width0.020em depth0em \kern-0.20em \rm C}
        }}

\newcommand{\wab}{W^{\alpha,\beta}}

\newcommand{\ba}{\begin{array}}
\newcommand{\ea}{\end{array}}
\newcommand{\be}{\begin{equation}}
\newcommand{\ee}{\end{equation}}
\newcommand{\bea}{\begin{eqnarray}}
\newcommand{\eea}{\end{eqnarray}}
\newcommand{\beaa}{\begin{eqnarray*}}
\newcommand{\eeaa}{\end{eqnarray*}}

\def\z{\zeta}

\font\tenmath=msbm10 \font\sevenmath=msbm7 \font\fivemath=msbm5
\newfam\mathfam \textfont\mathfam=\tenmath
\scriptfont\mathfam=\sevenmath \scriptscriptfont\mathfam=\fivemath

\def \={{\buildrel {\rm (law)} \over =}}

\def\qed{ \hfill \vrule width.25cm height.25cm depth0cm\smallskip}

\newcommand{\basa}{\begin{assumption}}
\newcommand{\easa}{\end{assumption}}

\newcommand{\bas}{\begin{assum}}
\newcommand{\eas}{\end{assum}}

\newcommand{\ignore}[1]{}
\begin{document}

\renewcommand{\thefootnote}{\fnsymbol{footnote}}

\title{{\bf Remarks on some linear fractional stochastic equations}}
\date{
\small {\bf Ivan Nourdin}\\
LPMA, Universit\'e Pierre et Marie Curie Paris 6,\\
Bo\^ite courrier 188, 4 Place Jussieu, 75252 Paris Cedex 5, France\\
{\tt nourdin@ccr.jussieu.fr}\\
$\mbox{ }$\\
$\mbox{ }$\\
\small {\bf Ciprian A. Tudor}\\
SAMOS/MATISSE, Universit\'e de Panth\'eon-Sorbonne Paris 1\\
90, rue de Tolbiac, 75634 Paris C\'edex 13, France\\
{\tt tudor@univ-paris1.fr}\\
}
 \maketitle

 \begin{abstract}
Using the multiple stochastic integrals we prove an existence and
uniqueness result for a linear stochastic equation driven by the
fractional Brownian motion with any Hurst parameter. We study  both
 the one parameter and two parameter cases. When the drift is zero,
 we show that in the one-parameter case the solution in an exponential, thus positive, function while in the
two-parameter settings the solution is negative on a non-negligible
set.
 \end{abstract}



\section{Introduction}
The significant amount of applications where the fractional Brownian
motion (fBm) is used led to the intensive development of the
stochastic calculus with respect to this process and its planar
version. The study of stochastic differential equations (SDEs) driven by
a fractional Brownian motion followed in a natural way. Let us
consider $(B^{\alpha }_{t}) _{t\in [0,T]} $ a fBm with Hurst
parameter $\alpha \in (0,1)$.
Essentially, one can consider the SDE
\begin{equation} \label{intro} dX_{t}=
\sigma (t,X_{t}) dB^{\alpha }_{t} + b(t,X_{t}) dt
\end{equation}
in two ways:
\begin{description}
\item{$\bullet$ } the pathwise (Stratonovich) type (that is, the
stochastic integral is considered in a pathwise sense);
\item{$\bullet$ } the divergence (Skorohod) type (that is, 
the stochastic integral is of divergence type).
\end{description}
The first type of equations, which includes the rough paths theory and the
stochastic calculus via regularization, can in general be solved
by now standard methods. We refer, among others, to  \cite{ALN,CQ,ER,FP,No,NoSi,NuRa,Za}. 
The second type
(Skorohod stochastic equations) is more difficult to be solved. Even
in the standard Brownian motion case (corresponding to $\alpha=1/2$),
we have an existence and uniqueness result only in two situations: 
\begin{description}
\item{$\bullet$} when $\sigma
(s,X_{s}) =\sigma(s) X_{s}$ with $\sigma(s)$ random: we then use an anticipating 
Girsanov transform, see \cite{buck},
\item{$\bullet$} when $\sigma(s,X_{s})
=\sigma(s)X_{s}$ and $b(s,X_{s}) =b(s)X_{s}$ with $\sigma,b$ two deterministic
functions: we can then use a method based on the Wiener-It\^o chaotic expansion.
\end{description}

This second approach will be considered in our paper. We will
consider the stochastic equation
\begin{equation}
\label{intro1} X_{t} =1+ \int_{[0,t]} aX_{s}\delta B^{\alpha }_{s}
+\int_{[0,t]}bX_{s}ds
\end{equation}
where $a,b$ are real numbers and the stochastic integral is
understood in the Skorohod sense. We first prove existence and
uniqueness results in the one-parameter case (that is when $t\in [0,T]$) and in
the two-parameter case (that is when $t\in [0,T]^{2}$ and with $B^{\alpha
}$ replaced by a fractional Brownian sheet $\wab$ with Hurst
parameters $\alpha,\beta$).

Of course, the fact that the above linear equation can be solved by
using Wiener-It\^o multiple integrals is not very surprising; it has
already used in \cite{PAT} for $\alpha >\frac{1}{2}$. Nevertheless,
we have to check some new technical aspects like: the proof of the
case $\alpha \in (0, \frac{1}{2})$ or the proof of the
two-parameter case for any Hurst parameters $\alpha $ and $\beta$.

More surprising is the behavior of the solution of (\ref{intro1})
when the drift $b$ is zero: in the one-parameter case, the solution
is an exponential, hence positive, function while
in the two-parameter case the solution is negative on a
non-negligible set. We also mention that, comparing to the standard
case when the Hurst parameters are $\frac{1}{2}$, new techniques like fractional Girsanov 
theorem and estimations of fractional norms of the kernels appearing in the
chaotic expression of the solution of (\ref{intro1}), are here needed.

We organized our paper as follows. Section 2 contains some
preliminaries on fractional Brownian motion and fractional Brownian
sheet. In Section 3 we study the existence, the uniqueness and the properties of the
solution of equation (\ref{intro1}) in both one-parameter and
two-parameter cases. Section 4 contains a technical proof.

\vskip0.5cm

\section{Preliminaries }
Consider  $(B^{\alpha }_{t}) _{t\in [0,T]} $ a fractional Brownian
motion (fBm) with Hurst parameter $\alpha \in (0,1)$   and let us
denote by $R^{\alpha }$ its covariance function
\begin{equation}
\label{cov1} R^{\alpha } (s,u)= \frac{1}{2}\left(
s^{2\alpha}+u^{2\alpha}-|s-u|^{2\alpha}\right)
\end{equation}
for every $s,u\in [0,T]$. It is well-known that $B^{\alpha } $
admits the Wiener integral representation
\begin{equation*}
B^{\alpha }_{t}= \int _{0}^{t} K^{\alpha } (t,s) dW_{s}
\end{equation*}
where $W$  denotes a standard Wiener process and
\begin{equation}
K^{\alpha }(t,s)=d_{\alpha}\left(  t-s\right)
^{\alpha-\frac{1}{2}}+s^{\alpha-\frac
{1}{2}}F_{1}\left(  \frac{t}{s}\right),   \label{for1}%
\end{equation}
$d_\alpha$ being a constant and $ F_1(z)=d_\alpha
\left(\frac{1}{2}-\alpha\right) \int_0^{z-1} \theta^{\alpha-3/2}
\left(1-(\theta+1)^{\alpha-1/2}\right)d\theta. $

 By
${\cal{H}}(\alpha )$ we will denote the canonical Hilbert space associated to
$B^{\alpha }$. That is, ${\cal{H}}(\alpha )$ is defined as the
closure of the linear space generated by the indicator functions $\{
1_{[0,t]}, t\in [0,T] \}$ with respect  to the scalar product
\begin{equation}
\label{prosca1} \langle 1_{[0,t]},1_{[0,u]}\rangle
_{{\cal{H}}(\alpha ) }= R^{\alpha } (t,u).
\end{equation}
The structure of ${\cal{H}}(\alpha )$  depends on
the values of the Hurst parameter $\alpha$. Let us recall the following
facts:
\begin{description}
\item{$\bullet$ }if $\alpha  \in (\frac{1}{2}, 1)$, then it
follows from \cite{PiTa} that the elements of ${\cal{H}}(\alpha )$
may be not functions but distributions of negative order. Thus it is
more convenient to work with subspaces of ${\cal{H}}(\alpha )$  that
are sets of functions. A such space is the set  $\vert
{\cal{H}}(\alpha )\vert $ of measurable functions
on $[0,T]$ such that
\begin{equation*}
\int _{0}^{T} \int _{0}^{T} \vert f(u) \vert \vert f(v)\vert \vert
u-v\vert ^{2\alpha -2} dudv <\infty
\end{equation*}
endowed with the scalar product
\begin{equation}
\label{hbar} \langle f, g\rangle _{ \vert {\cal{H}}(\alpha )\vert }=\alpha (2\alpha -1)\int _{0}^{T} \int _{0}^{T}f(u) f(v) \vert
u-v\vert ^{2\alpha -2} dudv.
\end{equation}

We have actually the inclusions
\begin{equation}
\label{inclu1} L^{2}([0,T]) \subset L^{\frac{1}{\alpha}}([0,T]) \subset
\vert {\cal{H}}(\alpha )\vert \subset {\cal{H}}(\alpha ).
\end{equation}

\item{$\bullet$ } if $\alpha \in (0,\frac{1}{2})$ then the Hilbert
space ${\cal{H}}(\alpha ) $ is a space of functions contained in
$L^{2} ([0,T])$. It contains the space of H\"older functions of
order $\alpha -\varepsilon $ with $\varepsilon >0$ and it can be
characterized by
\begin{equation}
\label{hk} {\cal{H}}(\alpha )= (K^{\ast })^{-1}(L^{2}([0,T]))
\end{equation}
where the operator $K^{\ast} $ is given by
\begin{equation}
\label{kstar} (K^{\ast} \varphi )(s) =K^\alpha(T,s)\varphi (s) + \int
_{s}^{T} (\varphi (r) -\varphi (s)) \frac{\partial K^\alpha}{\partial
r}(r,s)dr .
\end{equation}

\end{description}

A fBm being a Gaussian process, it is possible to construct multiple
Wiener-It\^o stochastic integrals with respect to it. We refer to
\cite{N} for general settings or to \cite{PAT} for the adaptation to
the fractional Brownian motion case. We only recall that the
multiple integral of order $n$ (denoted by $I_{n}$)  is an isometry
from ${\cal{U}}^{\otimes n} $ to $L^{2}(\Omega
)$ where ${\cal{U}}$ is the Hilbert space $\vert
{\cal{H}}(\alpha)\vert$  if $\alpha \in (\frac{1}{2}, 1)$ and the Hilbert space
${\cal{H}}(\alpha)$ if $\alpha \in  (0, \frac{1}{2})$.

We need to introduce the space $D^{ch}$ of stochastic processes that
can be expressed in terms of multiple stochastic integrals. That is,
we denote by $D^{ch}$ the set of
processes $u\in L^{2}(\Omega ; {\cal{U}} ) $
such that for every $t\in [0,T]$,
$$u_{t}= \sum _{n\geq 0} I_{n}(f_{n}(\cdot , t)) $$
where $f_{n} \in {\cal{U}}^{\otimes n+1} $ is
symmetric in the first $n$ variables and
\begin{equation}
\label{sko} \sum _{n\geq 1} (n+1)! \Vert f_{n} \Vert   ^{2} _{
{\cal{U}}^{\otimes n+1}  } <\infty.
\end{equation}
It follows from \cite{PAT} (for $\alpha >\frac{1}{2}$) or \cite{LN}
(for $\alpha <\frac{1}{2}$) that if $u\in D^{ch}$ then $u$ is
Skorohod integrable with respect to the fBm $B^{\alpha }$ and in
this case its Skorohod integral is
\begin{equation}
\label{delta2} \delta (u)=  \sum _{n\geq 0} I_{n+1} (\tilde{f}_{n})
\end{equation}
where $\tilde{f}_{n}$ means the symmetrization of $f_{n}$ with
respect to $n+1$ variables. Actually, in the case $\alpha
<\frac{1}{2}$ the expression (\ref{delta2}) corresponds to the
divergence integral in the extended sense.

\vskip0.5cm

Let us consider now the two-parameter case. Here, $\wab$ is a
fractional Brownian sheet with Hurst parameters $\alpha, \beta \in
(0,1)$. Recall that $\wab$ is defined as a centered Gaussian process
starting from $(0,0) $ with the covariance function

\begin{eqnarray*}
E\left(  W_{s,t}^{\alpha,\beta}W_{u,v}^{\alpha,\beta}\right)   &
=& R^{\alpha,\beta}(s,t,u,v)\\
&:=  &  \frac{1}{2}\left(
s^{2\alpha}+u^{2\alpha}-|s-u|^{2\alpha}\right) \frac{1}{2}\left(
t^{2\beta}+u^{2\beta}-|t-v|^{2\beta}\right)
\end{eqnarray*}
and it can be represented as
\[
W_{s,t}^{\alpha,\beta}=\int_{0}^{t}\int_{0}^{s}K^{\alpha}(t,u)K^{\beta
}(s,v)dW_{u,v}%
\]
where $(W_{u,v})_{u,v\in\lbrack0,T]}$ is a standard Brownian sheet
and $ K^{\alpha }$ is given by (\ref{for1}).  Denote by
\begin{equation*}
K^{\alpha , \beta }(t,s)=K^{\alpha}(t,u)K^{\beta }(s,v).
\end{equation*}
and let   ${\cal{H}}^{(2)}(\alpha , \beta ):= {\cal{H}}^{(2)} $ be  the
canonical Hilbert space of the fractional Brownian sheet $\wab$.
That is, ${\cal{H}}^{(2)}$ is defined as the closure of  the set of
indicator functions $\{ 1_{[0,t]\times [0,s]}, t,s \in [0,T]\} $
with respect to the scalar product
\begin{equation}
\label{prosca} \langle  1_{[0,t]\times [0,s]},  1_{[0,u]\times
[0,v]} \rangle _{{\cal{H}}^{(2)}}= R^{\alpha,\beta } (s,t,u,v)
\end{equation}
for every $t,s,u,v\in [0,T]$.

\smallskip

By the above considerations, we will have:
\begin{description}
\item{$\bullet$ }if $\alpha , \beta \in (\frac{1}{2}, 1)$, the elements of ${\cal{H}}^{(2)}$ may be
not functions but distributions. Thus it is more convenient to work
with subspaces of ${\cal{H}}^{(2)}$ that are sets of functions. We
have actually the inclusions
\begin{equation}
\label{inclu} L ^{2} ([0,T]^{2})
\subset \vert {\cal{H}}\vert ^{(2)} \subset {\cal{H}}^{(2)}
\end{equation}
where
$$\vert {\cal{H}}\vert ^{(2)} = \vert {\cal{H}}(\alpha )\vert \otimes
\vert {\cal{H}}(\beta )\vert$$ and $\vert {\cal{H}}(\alpha )\vert $
is defined by (\ref{hbar}).
\item{$\bullet$ }if $\alpha , \beta \in  (0, \frac{1}{2})$ then the
canonical space ${\cal{H}}^{(2)}$ is a space of functions that can be
written as
\begin{equation}
\label{hk2} {\cal{H}}^{(2)}= (K^{\ast, 2} )^{-1}\left(
L^{2}([0,T])\right) \subset L^{2}([0,T])
\end{equation}
where $K^{\ast ,2}$  is the product operator $K^{\ast }\otimes
K^{\ast} $ and $K^{\ast}$ is given by (\ref{kstar}).
\item{$\bullet$ } if $\alpha \in (\frac{1}{2}, 1)$ and $\beta \in
(0, \frac{1}{2})$, then $\vert {\cal{H}}\vert ^{(2)}$ is not a space of
functions and we will work with the subspace $\vert {\cal{H}}(\alpha
)\vert \otimes {\cal{H}}(\beta)$.
\end{description}
Let us denote by ${\cal{V}}$ the Hilbert space: $\vert
{\cal{H}}\vert ^{(2)}$  if $\alpha , \beta \in (\frac{1}{2}, 1)$;
${\cal{H}}^{(2)}$ if $\alpha , \beta \in  (0, \frac{1}{2})$ and $\vert
{\cal{H}}(\alpha )\vert \otimes {\cal{H}}(\beta)$ if  $\alpha \in
(\frac{1}{2}, 1)$ and $\beta \in (0, \frac{1}{2})$.

We can of course  consider multiple stochastic integrals with
respect to the Gaussian process $\wab$. Here the multiple integral
of order $n$, still denoted by $I_{n}$, will be a isometry from
${\cal{V}}^{\otimes n} $ to $L^{2}(\Omega)$.

\vskip0.5cm

\section{Linear stochastic equations with fractional Brownian motion and fractional Brownian sheet}

Let us consider the following stochastic integral equation
\begin{equation}\label{ec1}
X_{t}=1+ \int _{0}^{t} aX_{s}\delta B^{\alpha}_{s} + \int
_{0}^{t}bX_{s}ds, \hskip0.5cm t\in [0,T],
\end{equation}
where $a,b \in \mathbb{R}$ and the stochastic integral above is
considered in the Skorohod sense. We will first prove the existence
and the uniqueness of the solution of (\ref{ec1}), in the space
$D^{ch}$. For $\alpha
>\frac{1}{2}$ this has been proved in \cite{PAT}.

\begin{prop}
The equation (\ref{ec1}) admits an unique solution $X\in D^{ch}$
given by
\begin{equation} \label{sol1}
X_{t}= \sum _{n\geq 0} I_{n} (f_{n}(\cdot ,t))
\end{equation}
where the kernels $f_{n}$ are given by
$$f_{0}(t)= e^{bt}$$
and for every $n\geq 1$,
\begin{equation}\label{ker1}
f_{n}(t_{1}, \ldots , t_{n}, t) = \frac{a^{n}}{n!} e^{bt} 1_{[0,t]}
^{\otimes n} (t_{1}, \ldots , t_{n}).
\end{equation}
\end{prop}
{\bf Proof: } The expression (\ref{ker1}) of the kernels $f_{n}$
follows from Proposition 3.40 of \cite{PAT}. One can also compute it
easily by the recurrence relation
\begin{equation} \label{rec}
f_{0}(t)=e^{bt}, \hskip0.2cm f_{n}(t_{1}, \ldots , t_{n}, t) = a
\tilde{f}_{n-1} (t_{1}, \ldots , t_{n-1}, t_{n} ) 1_{[0,t]} (t_{n}), \hskip0.3cm \forall n\geq 1.
\end{equation}
  We only then
need to prove that $f_{n}\in \vert {\cal{H}}(\alpha )\vert ^{\otimes
n+1}$ (if $\alpha
> \frac{1}{2}$) and $f_{n}\in {\cal{H}}(\alpha ) ^{\otimes n+1}$
(if $\alpha < \frac{1}{2}$) and that the sum  (\ref{sko}) converges.

If $\alpha >\frac{1}{2}$, this follows easily from the inclusion
(\ref{inclu1}), since
\begin{equation*}\Vert \tilde{f}_{n} \Vert _{
\vert {\cal{H}}(\alpha )\vert ^{\otimes n+1}} \leq {\rm cst}\,\Vert \tilde{f}
_{n} \Vert _{L^{2}([0,T]^{n+1})}
\end{equation*}
and we can reduce to the classical situation ($\alpha =\frac{1}{2}$)
where the result is known.

If $\alpha <\frac{1}{2}$, then we need a new proof because the norm
${\cal{H}}(\alpha )$ is bigger than the norm $L^{2}$. Let us show
that the kernel $f_{n}$ given by (\ref{ker1}) (viewed as a function
of $n+1$ variables $t_{1}, \ldots t_{n}, t$) belongs to the space
${\cal{H}}(\alpha )^{\otimes n+1}. $ Here we can adapt an argument
used in \cite{TV2}. We will show that
\begin{equation*}
K^{\ast, n+1} f_{n} \in L^{2} ([0,T]^{n+1})
\end{equation*}
where $K^{\ast, n}  $ is the $n $ times tensor product of
$K^{\ast}$. It holds, by applying first the operator  $K^{\ast }$ to
the variables $t_{1}, \ldots , t_{n}$ and then to the variable $t$,
\begin{equation*}
K^{\ast, n+1} f_{n}= \frac{a^{n}}{n!} K^{\ast} e^{bt} \left(
K^{\ast , n} (1_{[0,t]}^{\otimes n}) \right)
\end{equation*}
and therefore, since
$$\Vert  K^{\ast , n} (1_{[0,t]}^{\otimes n}) \Vert ^{2}_{{\rm L}^2([0,T]^{n})}=
\Vert  K^{\ast} (1_{[0,t]}) \Vert ^{2n}_{{\rm L}^2([0,T]^{n})}
=t^{2n\alpha} $$ we get
\begin{eqnarray*}
\Vert K^{\ast, n+1} f_{n} \Vert  _{L^{2}
([0,T]^{n+1})} &=&\frac{|a|^{n}}{n!}\Vert     K^{\ast } (e^{bt}
t^{2n\alpha})\Vert _{L^{2}([0,T])} \\
&=& \frac{|a|^{n}}{n!} \Vert \sum _{k\geq 0 }
\frac{|b|^{k}}{k!}K^{\ast } \left( t^{k+2\alpha n}\right)\Vert 
_{L^{2}([0,T])}.
\end{eqnarray*}
 Since for every $k\geq 1$,
the function $t^{k+2\alpha n} $ is Lipschitz, then we have, using
(\ref{kstar})
\begin{eqnarray*}
\Vert  K^{\ast } \left( t^{k+2\alpha n}\right)\Vert _{L^{2}([0,T])}
&\leq & C(\alpha , T)(k+2\alpha n)  \left[ \int _{0}^{T} K^{\alpha}(T,t)^2 t^{4n\alpha
+ 2k} dt \right.\\
&&+\left. \int  _{0}^{T} \left( \int _{t}^{T} (r-t)
\frac{\partial K^\alpha}{\partial r} (r,t) dr \right) ^{2}dt\right]^{1/2} \\
&\leq & C(\alpha , T) (k+2\alpha n)T ^{k +\alpha (2n +1) }.
\end{eqnarray*}
This implies that
\begin{equation}
\label{a1} \Vert K^{\ast, n+1} f_{n} \Vert  _{L^{2} ([0,T]^{n+1})}
\leq {\rm cst} \frac{|a|^{n}}{n!}T^{\alpha (2n+1)}.
\end{equation}
The function $f_{n}$ being symmetric in the first $n$ variables, we have
\begin{equation*}
\tilde{f}_{n} (t_{1}, \ldots , t_{m+1}) =\frac{1}{m+1} \sum
_{i=0}^{m+1} f_{n} (t_{1}, \ldots , t_{m+1} ^{i} , \ldots , t_{m} )
\end{equation*}
where $t_{m+1} ^{i}$ means that $t_{m+1}$ is on the position $i$.
Clearly the bound (\ref{a1}) holds for $\tilde{f}_{n}$. By the above
estimate, it is not difficult to see that the sum (\ref{sko}) is convergent
because
\begin{equation*}
\sum _{n\geq 0} (n+1) ! \Vert \tilde{f}_{n}\Vert ^{2} _{
{\cal{H}}(\alpha ) ^{\otimes n+1} } \leq {\rm cst} \sum _{n\geq 0}
\frac{a^{2n}}{n! } T^{2\alpha  (2n+1)} <\infty. \end{equation*} The
uniqueness of the solution in $D^{ch}$ is obvious because,  if there
are two solutions, then the kernels of the chaotic expansion
verifies both the relation (\ref{rec}).  \qed

\vskip0.5cm

In the particular case when the drift $b$ is zero, we have the
following
\begin{corollary}\label{cor1}
The unique solution in $D^{ch}$  of the equation
\begin{equation}
\label{ec1''}X_{t}=1+ \int _{0}^{t} aX_{s}\delta B^{\alpha}_{s},
\hskip0.5cm t\in [0,T]
\end{equation}
is given by
\begin{equation} \label{exp}
X_{t}= \exp \left(  aB^{\alpha }_{t} -\frac{a^{2}}{2}t^{2\alpha }
\right).
\end{equation}
\end{corollary}
{\bf Proof: } This is of course a consequence of Proposition 1. But,
to compare to the two-parameter case, we prefer to show how the
formula (\ref{exp}) is obtained.  Let, for every $t\in [0,T]$
\begin{equation*}
X_{t}=\sum_{n\geq 0} I_{n}(f_{n}\left( \cdot , t\right) )
\end{equation*}
be the chaotic expression of $X$. Equation (\ref{ec1''}) can be
rewritten as
 \begin{equation}
 \label{ec2}
\sum _{n\geq 0} I_{n}\left( f_{n}(\cdot , t)\right) =1+a \sum
_{n\geq 0} I_{n+1}\left(   \widetilde{ f_{n}\left( \cdot, \star
\right) 1_{[0,t]}(\star) } \right)
 \end{equation}
 where $\cdot$ represents $n$ variables, $\star$ denotes one variable and
 $ \widetilde{ f_{n}\left( \cdot, \star
\right) 1_{[0,t]}(\star) }$ denotes the symmetrization of the
function $f_{n}\left( \cdot, \star \right) 1_{[0,t]}(\star)$ in
$n+1$ variables.

By identifying the corresponding Wiener chaos, we easily get
$$f_{0}(t)=1, \hskip0.3cm f_{1}(t_{1}, t) = a\,1_{[0,t]}(t_{1}) $$
and
$$f_{2}(t_{1}, t_{2}, t) = \frac{a^2}{2}\left(
1_{[0,t_{1}]}(t_{2})1_{[0,t]}(t_{1})+
1_{[0,t_{2}]}(t_{1})1_{[0,t]}(t_{2})\right) =\frac{a^2}{2}
1_{[0,t]}^{\otimes 2}(t_{1}, t_{2}).$$ By induction we will get for
every $n\geq 1$
\begin{equation}
\label{fn} f_{n}(t_{1}, t_{2}, \ldots , t_{n}) =\frac{a^n}{n!} \sum
_{i=0}^{n} 1_{[0, t_{i}]}^{\otimes n-1 }(\hat{t_{i}})1_{[0,t]}(t_i)
=\frac{a^n}{n!} 1_{[0,t]}^{\otimes n} (t_{1}, \ldots , t_{n})
\end{equation}
where by $\hat{t_{i}}$ we denoted the vector $(t_{1}, \ldots ,
t_{n})$ with $t_{i}$  missing. Therefore, we can express the
solution of (\ref{ec1''}) as
\begin{equation}
\label{xt} X_{t}= \sum _{n\geq 0} \frac{a^n}{n!}\,I_{n}\left(
1_{[0,t]}^{\otimes n} \right) = \exp \left( a\,B^{H}_{t}
-\frac{a^2}{2} \Vert  1_{[0,t]} \Vert ^{2} _{\cal{H}}\right)
\end{equation}
where for the last equality we refer e.g. to \cite{HDD}.\qed

\begin{remark}
{\rm In Skorohod setting, it is difficult, in general, to write an
Euler's type scheme associated to the equation $X_t=x_0+\int_0^t
\sigma(X_s)\delta B^\alpha_s$, even if $\alpha\ge 1/2$. Indeed, by using the integration by
parts for the Skorohod integral $\delta $ and the Malliavin
derivative $D$ (see \cite{N}) $\delta(Fu)=F\delta(u)-\langle DF,u
\rangle_{\cal{H}(\alpha)}$ and by assuming that we approximate
$X_{(k+1)/n}$ by $X_{k/n}+\int_{k/n}^{(k+1)/n} \sigma(X_{k/n})\delta
B^\alpha_s$ (as in the case $\alpha=1/2$), one obtains
$$
\widehat{X}^{(n)}_{(k+1)/n}=\widehat{X}^{(n)}_{k/n}
+\sigma( \widehat{X}^{(n)}_{k/n}) \left(B^\alpha_{(k+1)/n}-B^\alpha_{k/n}\right) - \sigma'(\widehat{X}^{(n)}_{k/n})
\langle D\widehat{X}^{(n)}_{k/n},1_{[k/n,(k+1)/n]}\rangle_{\mathcal{H}(\alpha)}.
$$
The problem is that the quantity $D\widehat{X}^{(n)}_{k/n}$ appears
and that it is difficult to compute it directly (without knowing the
solution).
Moreover, standard Euler
scheme do not apply here because the $L^{2}$-norm of the Skorohod
integral involves the first Malliavin derivative which involves the
second Malliavin derivative etc. and we cannot have closable
formulas. In the linear case, taking advantage from the
fact that we know explicitly the solution,  we can see what the
correct Euler scheme should be. Indeed, since we have $DX_{k/n}=aX_{k/n}{\bf
1}_{[0,k/n]}$ (see Corollary \ref{cor1} above), a natural Euler's type scheme associated to
(\ref{ec1''}) is
$$
\widehat{X}^{(n)}_{(k+1)/n}=\widehat{X}^{(n)}_{k/n}
+a \widehat{X}^{(n)}_{k/n} \left(B^\alpha_{(k+1)/n}-B^\alpha_{k/n}\right)
- \frac{a^2}{2}\widehat{X}^{(n)}_{k/n}\left[\left(\frac{k+1}{n}
\right)^{2H}-\left(\frac{k}{n}\right)^{2H}-\frac{1}{n^{2H}}\right].
$$
In fact, it is not very difficult to prove (using the same method as in the proof
of Proposition 6 in \cite{nourdincras}) that $(\widehat{X}^{(n)}_1)$
converges in ${\rm L}^2(\Omega)$ if and only if $\alpha\ge 1/2$ and that, in the case where $\alpha>1/2$, the limit is
$\exp \left(  aB^{\alpha }_{1} -\frac{a^{2}}{2}\right)$.
}
\end{remark}

As we have seen, the solution of (\ref{ec1''}) is an exponential,
hence positive, function. We will show that the situation is
different in the two-parameter case.

Before that, let us consider the equation corresponding to (\ref{ec1}) in the
two-parameter case
\begin{equation}\label{ec1'}
X_{z} = 1 + \int _{[0, z]}a X_{r}\delta \wab _{r}+ \int _{[0, z]}b
X_{r}dr
\end{equation}
where $z=(s,t)\in [0, T]^{2}$ and $\wab$ is a fractional Brownian
sheet with Hurst parameters $\alpha, \beta \in (0,1)$.

\vskip0.3cm

 We will denote now by $D^{ch, 2}$ the class of
functionals that can be represented as a serie of multiple
stochastic integrals with respect to $\wab$ (that is, $D^{ch,2}$ is
the two-parameter equivalent of $D^{ch}$). In the next proposition,
we show that (\ref{ec1'}) admits a unique solution in this space:

\begin{prop}
Let us denote by
$\mathcal{A}_n$ the set $\{(z_1,...,z_n)\in(\mathbb{R}^2)^n:\,\,\exists
\sigma\in\mathfrak{S}_n,\,z_{\sigma(1)}\le ...\le z_{\sigma(n)}\}$
(in the one-parameter case: $\mathcal{A}_n=\mathbb{R}^n$). If
$z\in\mathcal{A}_n$, we consider $\sigma=\sigma_z\in\mathfrak{S}_n$ such
that $z_{\sigma(1)}\le ...\le z_{\sigma(n)}$.

 The equation (\ref{ec1'}) admits an unique solution $X\in D ^{ch,2}$ given by $X_{z}=\sum_{n\geq 0} I_{n}
( f_{n}(\cdot , z)) $ where
\begin{equation}\label{fn'}
\begin{array}{llll}
f_n(z_1,...,z_n,z)&=&\frac{a^n}{n!}
h_0(b(s-s_{\sigma_z(n)})(t-t_{\sigma_z(n)}))
\times {\bf 1}_{\mathcal{A}_n}(z_1,...,z_n){\bf 1}^{\otimes n}_{[0,z]}(z_1,...,z_n)\\
&&\times\prod_{1\le j\le n}
h_0(b(s_{\sigma_z(j)}-s_{\sigma_z(j-1)})(t_{\sigma_z(j)}-t_{\sigma_z(j-1)}))
\end{array}
\end{equation}
with $z=(s,t)$, $z_i=(s_i,t_i)$ and $h_0(x)=\sum_{n=0}^\infty
\frac{x^n}{(n!)^2}$. We also used the convention that $\sigma_z(0)=0$ and $z_0=(0,0)$.\\
\end{prop}
{\bf Proof: } We only prove the algebraic part (\ref{fn'}) of the Proposition.
Indeed, the fact that the kernels
$f_{n}$ belongs to ${\cal{V}}^{\otimes n+1} $ did not present new
difficulties with respect to the proofs of Propositions 1 and 3.
Thus, we return to these proofs for this point.
Let us write
$$X_{z}= \sum _{n\geq 0} I_{n}\left( f_{n} (\cdot, z) \right).  $$
Here, $I_{n}$ is the $n$-order Wiener-It\^o multiple integral with
respect to the fractional Brownian sheet $W^{\alpha,\beta}$ and $f_{n}\in L^{2}\left(
[0,T]^{2n}\right)$.
From (\ref{ec1'}) we have that $f_{0} (z) =h_0(bst)$ and
for $n\geq 1$,
$$
f_n(z_1,...,z_n,z)=a \widetilde{f_{n-1}(z_1,...,z_n) 1_{[0,z]}(z_n)}
+ b \int_{[0,z]}f_n(z_1,...,z_n,r)dr.\\
$$
Let  $n=1$. We therefore have
$$
f_1(z_1,z)=a\,h_0(bs_1t_1){\bf 1}_{[0,z]}(z_1)  +
b\int_{[0,z]}f_1(z_1,r)dr
$$
and
$$
f_1(z_1,z)=a\,h_0(bs_1t_1)\,h_0\left(b(s-s_1)(t-t_1)\right) {\bf
1}_{[0,z]}(z_1)
$$
hence (\ref{fn'}) is satisfied. If $n=2$ it holds that
$$
\widetilde{f_{2}(z_1,z_2,z_3){\bf 1}_{[0,z]}(z_3)}=\frac{1}{2}
\left( a\,h_0(bs_1t_1)\,h_0\left(b(s_2-s_1)(t_2-t_1)\right) {\bf
1}_{0\le z_1\le z_2\le z} \right.
$$
$$
\left. + a\,h_0(bs_2t_2)\,h_0\left(b(s_1-s_2)(t_1-t_2)\right) {\bf
1}_{0\le z_2\le z_1\le z} \right).
$$
Since
$$
f_2(z_1,z_2,z)=a\widetilde{f_{1}(z_1,z_2){\bf 1}_{[0,z]}(z_2)}   +
b\int_{[0,z]}f_1(z_1,r)dr
$$
we deduce that
$$
f_2(z_1,z_2,z)=\frac{a^2}{2} \left(
h_0(bs_1t_1)\,h_0\left(b(s_2-s_1)(t_2-t_1)\right)\,h_0\left(b(s-s_2)(t-t_2)\right)
{\bf 1}_{0\le z_1\le z_2\le z}\right.
$$
$$
\left.+
h_0(bs_2t_2)\,h_0\left(b(s_1-s_2)(t_1-t_2)\right)\,h_0\left(b(s-s_1)(t-t_1)\right)
{\bf 1}_{0\le z_2\le z_1\le z} \right)
$$
and again  (\ref{fn'}) is  verified. The above computations can be
easily extended to an induction argument. \qed

\vskip0.5cm

Let us now discuss the case $b=0$:
\begin{prop}
The equation
\begin{equation}\label{eq1}
X_{z} = 1 + \int _{[0, z]}a X_{r}\delta \wab _{r}, \hskip0.5cm z\in
[0,T]^{2}
\end{equation}
admits  an unique solution $X\in D ^{ch,2}$ given by $X_{z}=\sum_{n\geq 0} I_{n}
( f_{n}(\cdot , z)) $ where
\begin{equation}
\label{fn2} f_ {n} \left( \rho _{1}, \ldots , \rho _{n}, z\right)=\frac{a^n}{n!}\sum _{i=1}^{n} 1_{[0, \rho _{i}]}^{\otimes n-1 }
(\hat{\rho _{i}})1_{[0,z]}(\rho _{i}).
\end{equation}
\end{prop}
{\bf Proof: }Let us write
$$X_{z}= \sum _{n\geq 0} I_{n}\left( f_{n} (\cdot, z) \right).$$
From the equivalent of relation (\ref{ec2})
in the two-parameter case, we obtain
$$f_{0}(z)=1, \hskip0.3cm f_{1}(\rho_1, z) = a\,1_{[0,z]}(\rho_1 )$$
and in general relation (\ref{fn2}) holds.  Since $\mathcal{A}_n\not =(\mathbb{R}^2)^n$ (recall
that $\mathcal{A}_n$ is defined in Proposition 2), note that this last
expression is not equal to $\frac{a^m}{m!}1_{[0,z]}^{\otimes m}
\left( \rho _{1}, \ldots \rho _{n}\right) $ as in the one-parameter
case (see Corollary 1).

Let us now prove that the kernel $f_{n}$ belongs to the space
${\cal{V}}^{\otimes n+1 }$. When the Hurst parameters $\alpha $ and
$\beta $ are bigger than $\frac{1}{2}$, then we can use
(\ref{inclu}) and then refer to the standard case of the Brownian
sheet. We will thus only discuss the case $\alpha , \beta
<\frac{1}{2}$; the case $\alpha >\frac{1}{2} $ and $\beta
<\frac{1}{2}$ will be a mixture  of the other  two cases. We use the
induction. We will illustrate first the case $n=2$. We check that
$1_{[0, z]}(z_{2}) 1_{[0, z_{2}]}(z_{1}) $ belongs to
${{\cal{H}}^{(2)}}^{\otimes 3}$. This actually reduces to proving
that
$$1_{[0,t]} (t_{2})1_{[0, t_{2}] } (t_{1}) \in  {\cal{H}}(\alpha ) ^{\otimes 3}.$$
Let us apply the operator $K^{\ast, 3}$ in three steps: first to the
variable $t_{1}$, then to the variable $t$ and then to $t_{2}$. It
holds that \begin{eqnarray*}
 \Vert K^{\ast ,3}\left(  1_{[0,t]}
(t_{2})1_{[0, t_{2}] } (t_{1})\right) \Vert^2  _{L^{2} ([0,T]^{3})}
&=& \Vert K^{\ast ,2}\left( t_{2} ^{2\alpha }1_{[0,t]} (t_{2})
\right)
\Vert^2  _{L^{2} ([0,T]^{2})}\\
&=& \Vert K^{\ast ,1} \left( t_{2} ^{2\alpha } (T-t_{2} ) ^{2\alpha
}\right)  \Vert^2 _{L^{2} ([0,T])}
\end{eqnarray*}
and to conclude we refer to Proposition 3.6 in \cite{CN}: it is a
straightforward consequence of Lemma 4.3 in \cite{CN} that
$(T-t_{2}) ^{2\alpha } (B^\alpha)^{2}$ belongs to the extended
domain of the divergence and therefore its expectation is in
${\cal{H}}(\alpha )$.

  We will show now  that the kernel $ 1_{[0,
\rho _{i}]}^{\otimes n-1 } (\hat{\rho _{i}})1_{[0,z]}(\rho _{i})$
has a finite norm in ${{\cal{H}}^{(2)}}^{\otimes n+1}$ by assuming
that the result is true for $n $ variables. In suffices to check
that the function of $n+1$ (real) variables
$$1_{[0,t]}(t_{n})1_{[0,t_{n}]}(t_{n+1})\ldots 1_{[0,
t_{2}]}(t_{1})$$ belongs to ${\cal{H}}(\alpha )^{\otimes n+1}$ or,
equivalently, the operator $K^{\ast, n+1} $ applied to the above
function is in $L^{2}([0,T]^{n+1})$. By applying first the operator
$K^{\ast}$ to the variable $t_{1}$ it holds that
\begin{eqnarray*}
&& \Vert K^{\ast , n+1} \left(
1_{[0,t]}(t_{n})1_{[0,t_{n}]}(t_{n+1})\ldots 1_{[0, t_{2}]}(t_{1})
\right) \Vert^{2} _{L^{2}([0,T]^{n+1})}\\
 &=&\Vert K^{\ast , n}
\left( 1_{[0,t]}(t_{n})1_{[0,t_{n}]}(t_{n+1})\ldots 1_{[0,
t_{3}]}(t_{2})t_{2} ^{2\alpha }  \right) \Vert^{2}
_{L^{2}([0,T]^{n})}\\
&=& \Vert K^{\ast }\left( t_{2}^{2\alpha } g(t_{2})\right) \Vert
^{2} _{L^{2}([0,T])}
\end{eqnarray*}
where the function $t_{2} \to g(t_{2}) :=\Vert K^{\ast , n-1} \left(
1_{[0,t]}(t_{n})1_{[0,t_{n}]}(t_{n+1})\ldots 1_{[0, t_{3}]}(t_{2})
\right) \Vert^{2} _{L^{2}([0,T]^{n-1})}$ belongs to
${\cal{H}}(\alpha )$ by the induction hypothesis. Now, we refer to
the proof of Proposition 3.6 in \cite{CN} for the fact  that
$g(\cdot ) E\left( B. ^{2}\right) $ has a finite norm in
${\cal{H}}(\alpha )$.

It can  actually be proved as above that $$\Vert 1_{[0, \rho
_{i}]}^{\otimes n-1 } (\hat{\rho _{i}})1_{[0,z]}(\rho _{i})\Vert _{
{{\cal{H}}^{(2)}}^{\otimes n+1} }\leq \frac{C^{n}}{n!}  $$ for every
$n$ where $C$ is a positive constant. Now  we can finish as in proof
of Proposition 1.\qed

\vskip0.5cm

We will need the following Girsanov theorem. Its proof will be given
in the Appendix.
\begin{lemma}\label{girs}
For any $\varepsilon >0$, the process
\begin{equation}
\label{we} W^{\alpha , \beta, \varepsilon }_{s,t} = \wab _{s,t}
-\frac{st}{\varepsilon }
\end{equation}
has the same law as a fractional Brownian sheet with parameters
$\alpha ,\beta $ under the new probability $P_{\varepsilon }$ given
by
\begin{equation}
\label{pe} \frac{dP^{\varepsilon }}{dP} =\exp \left(
\frac{1}{\varepsilon }\wab _{T,T} -\frac{1}{2\varepsilon^{ 2}}\int
_{[0,T]^{2}}\left( K_{\alpha, \beta } ^{-1}\left( F(\cdot )\right)
(\rho ) \right) ^{2} d\rho \right)
\end{equation}
where $F(t,s)=ts$ and $K_{\alpha , \beta }$ is the operator
associated to the kernel of the $\wab$.
\end{lemma}

\vskip0.5cm

The solution of the equation (\ref{eq1}) has actually a different
behavior comparing to the one-parameter case (Corollary \ref{cor1}). We prove
actually below that the solution of (\ref{eq1}) is almost surely
negative on a non-negligible set. Note that the same problem has
been studied in the case of the standard Brownian sheet in
\cite{N87}.

\begin{prop}\label{prop1}
Let $X$ be the unique solution to (\ref{eq1}) in the space $D^{ch,2}$. Then
there exists an open set $\Delta \subset [0,T]^{2}$ such that
\begin{equation}
\label{neg} P\{ X_{z}<0  \mbox{ for all } z\in \Delta \} >0.
\end{equation}
\end{prop}
{\bf Proof: } Note that the deterministic equation
\begin{equation}
\label{ecdet} g(s,t)= 1+ \int_{0}^{s}\int_{0}^{t} ag(u,v) dudv
\end{equation}
admits the unique solution $g(s,t)=h_{0}(ast) $ with
$h_{0}(x)=\sum_{n\geq 0} \frac{x^{n}}{(n!)^{2}}$ and that the
function $h_{0}$ satisfy the property: there exists an open set
$I=(-\beta , -\alpha )$ such that $h_{0} (x) <-\delta <0$ for any
$x\in I$ (see \cite{N87}, page 231).

Suppose $a>0$, fix $N>0$ and define the open set
\begin{equation*}
\label{delta} \Delta =\{ (s,t), \alpha < ast< \beta , 0< s,t<N\}.
\end{equation*}
For every $\varepsilon >0$, consider
\begin{equation*}
\label{ec3} X_{z}^{\varepsilon } =1 + \int _{[0,z]} a\varepsilon
X^{\varepsilon }_{r} \delta\wab _{r}.
\end{equation*}
Thanks to Corollary \ref{cor1}, we know that the solution $X^{\varepsilon}$ of (\ref{ec3}) is
given by
\begin{equation*}
\label{ec4} X_{z}^{\varepsilon }=\sum_{n\geq 0} \varepsilon
^{n} I_{n} \left( f_{n}(\cdot ,z)\right)
\end{equation*}
where the kernels $f_{n}$ are given by $ f_ {n} \left( \rho _{1},
\ldots , \rho _{n}, z\right)=\frac{a^n}{n!}\sum _{i=1}^{n} 1_{[0,
\rho _{i}]}^{\otimes n-1 } (\hat{\rho _{i}})1_{[0,z]}(\rho _{i}). $
 Let us consider the equation
\begin{equation}\label{ec5}
Y^{\varepsilon}_{z}= 1 + \int _{[0,z]} a\varepsilon Y^{\varepsilon
}_{r}dW^{\alpha , \beta , \varepsilon }_{r}= 1 +\int _{[0,z]}
a\varepsilon Y^{\varepsilon }_{r}dW^{\alpha , \beta }_{r}-\int
_{[0,z]} a Y^{\varepsilon }_{r}dr
\end{equation}
and recall that, by Lemma \ref{girs}, $W^{\alpha , \beta , \varepsilon
}$ is a fractional Brownian sheet under $P_{\varepsilon}$.  Now, we
observe that
\begin{equation}\label{Kfini}
K=\sup _{\varepsilon >0}\sup _{z}E\left| Y^{\varepsilon }_{z}\right|
^{2}<\infty.
\end{equation}
In fact, to show that (\ref{Kfini}) holds is not difficult because
it follows from Proposition 2 that the kernel of order $n$ appearing
in the chaotic expression of the solution of (\ref{ec5}) are of the
form $\varepsilon ^{n}$ multiplied to the kernel of order $n$ of the
solution of (\ref{ec5}) with $\varepsilon =1$. Then, $K\le \sup
_{z}E\left| Y^{1 }_{z}\right| ^{2}<\infty$.

Now, by (\ref{ecdet}) and (\ref{ec5}) we have, if $z=(t,s)$,
\begin{equation*}
 Y^{\varepsilon }_{s,t} -h_{0}(-ast) = -a\int _{0}^{s} \int
_{0}^{t} \left(Y^{\varepsilon }_{u,v }-h_{0}(-auv)\right) dvdu +
a\varepsilon \int _{0}^{s}\int _{0}^{t}Y^{\varepsilon }_{u,v
}dW^{\alpha , \beta }_{u,v}
\end{equation*}
and using the bound (\ref{Kfini}) and the Gronwall Lemma in the
plane we obtain
$$E \left[\sup_{s,t\in [0,T]}\left|  Y^{\varepsilon }_{s,t} -h_{0}(-ast)
\right| ^{2}\right] \to _{\varepsilon \to 0} 0 .$$
Since $h_{0}(-ast) <-\delta $ for $(s,t)\in \Delta $, it follows
that
\begin{equation*}
P\left(   Y^{\varepsilon} _{z}<0, \forall z\in \Delta \right) \to _{\varepsilon \to 0}1.
\end{equation*}
Thus, for every $\varepsilon>0$ small enough
$$P\left(   Y^{\varepsilon} _{z}<0, \forall z\in \Delta \right) >0 $$
and
$$P_{\varepsilon }\left(   Y^{\varepsilon} _{z}<0, \forall z\in \Delta \right)
=P\left(   X^{\varepsilon} _{z}<0, \forall z\in \Delta \right)
>0.
$$
Since $\wab _{c_{1}s,c_{2}t}$ has the same law as
$c_{1}^{\alpha }c_{2}^{\beta} \wab_{s,t} $ as process, we get that
$X^{\varepsilon }_{s,t} $ is equal in law to $X_{\varepsilon
^{2\alpha} s, \varepsilon ^{2\beta }t}$. So, for $\varepsilon>0$ small enough,
$$P\left( X_{\varepsilon
^{2\alpha }s , \varepsilon ^{2\beta }t} <0, \forall z\in \Delta
\right) >0$$ and the conclusion follows.\qed

\section{Appendix}

{\bf Proof of Lemma \ref{girs}: }The conclusion will follow from the
Girsanov theorem for the fractional Brownian sheet (see Theorem 3 in
\cite{ENO}) if we show that the functions $F(s,t) =st$ belongs to
the space $I^{\alpha + \frac{1}{2}, \beta + \frac{1}{2}}\left(
L^{2}([0,T]^{2})\right) $ or equivalently,
\begin{equation*}
K_{\alpha, \beta } ^{-1}\left( F(\cdot )\right)\in L^{2}
([0,T]^{2}).
\end{equation*}
To show this, we will need the expression of its inverse operator in
terms of fractional integrals and derivatives (see e.g. \cite{ENO})
\begin{equation}
\label{k1} K_{\alpha ,\beta }^{-1} h(t,s)= t^{\alpha
-\frac{1}{2}}s^{\beta -\frac{1}{2}}I^{\frac{1}{2}-\alpha ,
\frac{1}{2}-\beta } \left( t^{\frac{1}{2}-\alpha
}s^{\frac{1}{2}-\beta }\frac{\partial ^{2}h}{\partial s\partial
t}\right), \hskip0.5cm \alpha ,\beta < \frac{1}{2}
\end{equation}
and
\begin{equation} \label{k2} K_{\alpha ,\beta }^{-1} h(t,s)=t^{\alpha
-\frac{1}{2}}s^{\beta -\frac{1}{2}}D^{\alpha-\frac{1}{2} ,\beta
-\frac{1}{2}}\left(t^{\frac{1}{2}-\alpha }s^{\frac{1}{2}-\beta
}\frac{\partial ^{2}h}{\partial s\partial t}\right),\hskip0.5cm
\alpha ,\beta > \frac{1}{2}.
\end{equation}
Here,
$$
I^{\alpha,\beta}f(x,y)=\frac{1}{\Gamma(\alpha)\Gamma(\beta)}\int_0^x\int_0^y
(x-u)^{\alpha-1}(y-v)^{\beta-1} f(u,v)dudv
$$
and
$$
D^{\alpha,\beta}f(x,y)=\frac{1}{\Gamma(1-\alpha)\Gamma(1-\beta)}
\frac{\partial^2}{\partial x\partial y} \int_0^x\int_0^y
\frac{f(u,v)}{(x-u)^{\alpha}(y-v)^{\beta}}dudv,
$$
with $\Gamma$ the Euler function.

For $\alpha , \beta  \in (0, \frac{1}{2}) $ we have

\begin{eqnarray*}
K_{\alpha, \beta } ^{-1} F (t,s)&=& t^{ \alpha-\frac{1}{2}
}s^{\beta-\frac{1}{2}}I^{\alpha -\frac{1}{2} ,\beta
-\frac{1}{2}}\left( t^{\frac{1}{2}-\alpha
}s^{\frac{1}{2}-\beta}\right)\\
&=& \frac{1}{\Gamma(\frac{1}{2}-\alpha )\Gamma(\frac{1}{2}-\beta )}
 t^{\alpha -\frac{1}{2} }s^{\beta
 -\frac{1}{2}}\int_{0}^{t}\int_{0}^{s} (t-u) ^{-\frac{1}{2}-\alpha
 } u^{\frac{1}{2}-\alpha }(s-v)^{-\frac{1}{2}-\beta
 }v^{\frac{1}{2}-\beta }dvdu
 \end{eqnarray*}
 and this belongs to $L^{2}
([0,T]^{2}).$

\smallskip

If $\alpha , \beta \in (\frac{1}{2}, 1)$ then by (\ref{k2}) we can
write
\begin{eqnarray*}
&&K_{\alpha, \beta } ^{-1}F (t,s)\\
&=& t^{\alpha -\frac{1}{2}}s^{\beta -\frac{1}{2}}D^{\alpha
-\frac{1}{2}, \beta -\frac{1}{2}}\left( t^{\frac{1}{2}-\alpha
}s^{\frac{1}{2}-\beta } \right) (t,s)\\
&=& t^{\alpha -\frac{1}{2}}s^{\beta
-\frac{1}{2}}\frac{1}{\Gamma(\frac{3}{2}-\alpha )\Gamma
(\frac{3}{2}-\beta )} \left[ \frac{ t^{\frac{1}{2}-\alpha
}s^{\frac{1}{2}-\beta  } }{ t^{\alpha
-\frac{1}{2}}s^{\beta -\frac{1}{2}}} \right. \\
&&\left. +  \frac{\alpha -\frac{1}{2}}{s^{\beta -\frac{1}{2}}}\int
_{0}^{t} \frac{ t^{\frac{1}{2}-\alpha }s^{\frac{1}{2}-\beta }-
u^{\frac{1}{2}-\alpha }s^{\frac{1}{2}-\beta
}}{(t-u)^{\alpha +\frac{1}{2}}}du    \right. \\
&&\left. + \frac{\beta -\frac{1}{2}}{t^{\alpha
-\frac{1}{2}}}\int_{0}^{s}  \frac{ t^{\frac{1}{2}-\alpha
}s^{\frac{1}{2}-\beta }-t^{\frac{1}{2}-\alpha }v^{\frac{1}{2}-\beta
}} {(s-v)^{\beta +\frac{1}{2}}}dv
\right. \\
&&\left. + (\alpha -\frac{1}{2})(\beta -\frac{1}{2})
\int_{0}^{t}\int_{0}^{s} \frac{ t^{\frac{1}{2}-\alpha
}s^{\frac{1}{2}-\beta }- u^{\frac{1}{2}-\alpha }s^{\frac{1}{2}-\beta
} - t^{\frac{1}{2}-\alpha }v^{\frac{1}{2}-\beta
}+u^{\frac{1}{2}-\alpha }v^{\frac{1}{2}-\beta }}{(t-u)^{\alpha
+\frac{1}{2}} (s-v)^{\beta +\frac{1}{2}}} dvdu \right].
\end{eqnarray*}
Since $$\int_{0}^{t}\frac{ t^{\frac{1}{2}-\alpha
}-u^{\frac{1}{2}-\alpha }}{(t-u)^{\alpha +\frac{1}{2}}}du =c(\alpha
) t^{1-2\alpha },$$ it is not difficult to see that the above
function is in  $L^{2} ([0,T]^{2}).$

If $\alpha \in (0, \frac{1}{2})$ and $ \beta \in (\frac{1}{2}, 1)$,
then we have \begin{eqnarray*}
 K_{\alpha ,
\beta } ^{-1}F(t,s) &=& C(\alpha , \beta ) t^{\alpha
-\frac{1}{2}}\int _{0}^{t}(t-u) ^{-\frac{1}{2}-\alpha
 } u^{\frac{1}{2}-\alpha }du \\
 &&\times s^{\frac{1}{2}-\beta } + (\beta -\frac{1}{2})\int_{0}^{s}  \frac{ t^{\frac{1}{2}-\alpha
}s^{\frac{1}{2}-\beta }-t^{\frac{1}{2}-\alpha }v^{\frac{1}{2}-\beta
}} {(s-v)^{\beta +\frac{1}{2}}}dv
\end{eqnarray*}
and the conclusion is clearly a consequence of the above two cases.
The proof of Lemma is done.\qed

\end{document}